# EXPLICIT CALCULATION OF SINGULAR INTEGRALS OF TENSORIAL POLYADIC KERNELS

By

M. PERRIN (*Université de Bordeaux, CNRS, LOMA, UMR 5798, 351 Cours de la libération, 33400 Talence, France*)

AND

F. GRUY ( *Mines Saint-Etienne, Univ. Lyon, CNRS, UMR 5307 LGF, Centre SPIN 42023 Saint-Etienne France* )

**Abstract.** The Riesz transform of $u : \mathcal{S}(\mathbb{R}^n) \to \mathcal{S}'(\mathbb{R}^n)$ is defined as a convolution by a singular kernel, and can be conveniently expressed using the Fourier Transform and a simple multiplier. We extend this analysis to higher order Riesz transforms, i.e. some type of singular integrals that contain tensorial polyadic kernels and define an integral transform for functions $\mathcal{S}(\mathbb{R}^n) \to \mathcal{S}'(\mathbb{R}^{n \times n \times \ldots n})$. We show that the transformed kernel is also a polyadic tensor, and propose a general method to compute explicitely the Fourier mutliplier. Analytical results are given, as well as a recursive algorithm, to compute the coefficients of the transformed kernel. We compare the result to direct numerical evaluation, and discuss the case $n = 2$, with application to image analysis.

**1. Introduction.** The Riesz transform is a particular case of singular integral transform that attracted a great interest in fundamental and applied mathematics (especially for image analysis and processing) [1–4].

It is defined by the convolution $u \stackrel{R_i}{\mapsto} u \star W_i$, where $W_i$ is the distribution :

$$\langle W_i, \varphi \rangle = \frac{\Gamma\left(\frac{n+1}{2}\right)}{(\pi)^{\frac{n+1}{2}}} \lim_{\varepsilon \to 0} \int_{\mathbb{R}^n \setminus B(0,\varepsilon)} \frac{x_i/\|x\|}{\|x\|^n} \varphi(x)\, dx, \qquad (1.1)$$

with $\varphi \in \mathcal{S}(\mathbb{R}^n)$ and $x = (x_1, \ldots, x_i, \ldots, x_n) \in \mathbb{R}^n$.[1]

The convolution that defines $R_i$, see [5, Eq. (5.1.39)], is expressed simply using the Fourier transform (FT), denoted $\mathcal{F}$, as :

$$\mathcal{F}R_i(u)(\xi) = \left(-i\frac{\xi_i}{\|\xi\|}\mathcal{F}(u)(\xi)\right), \; \textit{for } u \in \mathcal{S}(\mathbb{R}^n). \qquad (1.2)$$

---

*E-mail address*: `mathias.perrin@u-bordeaux.fr`

[1]In the following, we shall also denote $\theta = (\theta_1, \theta_2, \ldots \theta_n)$ and $\xi = (\xi_1, \xi_2, \ldots, \xi_n)$, two elements of $\mathbb{R}^n$.

Therefore, in the Fourier domain, the Riesz transform is basically operated by the multiplier $-i\xi/\|\xi\|$ (using vectorial notation [4] to encompass all the components $x_i$).

In this article, we will generalize this approach, considering a class of singular integral transforms, whose singular kernel also varies as $x \mapsto 1/\|x\|^n$, and that can be found in electromagnetic or hydrodynamic problems. It is obtained replacing $(W_i)_{i \in [1..n]}$ by $W_f$,

$$\langle W_f, \varphi \rangle = \lim_{\varepsilon \to 0} \int_{\mathbb{R}^n \setminus B(0,\varepsilon)} \frac{f(x/\|x\|)}{\|x\|^n} \varphi(x)\, dx, \tag{1.3}$$

where $f$ is an integrable function $\mathbb{S}^{n-1} \to \mathbb{R}^{\underbrace{n \times n \times \cdots \times n}_{t \text{ times}}}$, whose image is an order $t$ and rank 1 tensor.

$\langle W_f, \varphi \rangle$ is a short way of writting $\langle W_{f_{i_1,i_2,\ldots,i_t}}, \varphi \rangle_{(i_1,i_2,\ldots,i_t) \in [1..n]^t}$, where $f_{i_1,i_2,\ldots,i_t} : \mathbb{S}^{n-1} \to \mathbb{R}$ is one component of the image tensor of $f$.

The components of the polyadic tensor we will consider are written [6] :

$$f_{i_1,i_2,\ldots,i_t}(\theta) = \prod_{q=1}^{t} \theta_{i_q} \tag{1.4}$$

or in short, $f(\theta) = \theta \circ \theta \circ \cdots \circ \theta$.

If we assume that

$$\int_{\mathbb{S}^{n-1}} f(\theta)\, d\theta = 0, \tag{1.5}$$

as a necessary condition, then, it is known [5], that the Fourier Transform (FT) of $W_f$ is defined on $\mathbb{R}^n$ by

$$\widehat{W_f}(\xi) = \int_{\mathbb{S}^{n-1}} f(\theta) \left[ -\ln(|\xi \bullet \theta|) + (-i\pi/2)\, sgn(\xi \bullet \theta) \right] d\theta, \tag{1.6}$$

where $\bullet$ is the canonical scalar product.

The purpose of this article is to compute explicitely $\widehat{W_f}(\xi)$, using integrals over $\mathbb{S}^{n-1}$, so as to generalise Eq. (1.2) to all the possible polyadic kernels with a $x \mapsto 1/\|x\|^n$ singularity, and thereby extend the notion of Riesz transform to higher order kernels.

It is convenient to rewrite the two integrals that appear in Eq. (1.6) as :

$$T(\xi) = \int_{\mathbb{S}^{n-1}} f(\theta)\, g(\xi \bullet \theta)\, d\theta \tag{1.7}$$

with $g(x) = sgn(x)$ or $g(x) = -\ln(|x|)$. In the following, we shall deal with $T$ instead of $\widehat{W_f}$.

In Sec. (2), we will formulate the problem, introducing a change of basis by Schmidt orthonormalization, to express $T$ in a new basis. In Sec. (3), we will make the computation, for which the key point is to factorize $T$ in the new basis. We shall show that, to each component of the $f(\theta)$ tensor, corresponds a component of the $T(\xi)$ tensor, which is also a sum of rank one tensors of order smaller or equal to $t$.

In Sec. (4), we present an algorithm to make the calculation in the general case, and compare numerical and analytical results. Finally, in Sec. (5), we will discuss as an example the case $n = 2$, and show that this singular transform can be used to process images. Perspectives will be drawn in the conclusion.

**2. Notations, change of basis.** It is simple to show [5], that $\widehat{W_f}(\xi) = \widehat{W_f}(\xi/|\xi|)$. Therefore, in the following, we shall assume that $|\xi| = 1$.

As the vector $\xi \in \mathbb{R}^n$ plays a particular role, we consider a change to a new basis, in which the first vector is $\xi$, and the $n-1$ remaining basis vectors are orthogonal to each other, and with $\xi$. The matrix for change of basis is denoted $R_\xi$ (we shall drop the index $\xi$ in the following). The detail of calculation is given in Sec. (7).

In the canonical basis, $T$ is defined by Eq. (1.7), that writes :

$$T = \int_{\mathbb{S}^{n-1}} \underbrace{\theta \circ \theta \circ \cdots \circ \theta}_{t \text{ terms}} \, g(\xi \bullet \theta) \, d\theta, \tag{2.1}$$

To simplify the integrand, it is interesting to represent the tensor $T$ in the new basis, where its expression becomes:

$$T' = \int_{\mathbb{S}^{n-1}} \underbrace{\theta' \circ \theta' \circ \cdots \circ \theta'}_{t \text{ terms}} g\left(\theta'_1\right) d\theta', \tag{2.2}$$

In practice, $\theta$ and $\theta'$ are expressed in spherical coordinates [5, Appendix D.].

Let us consider a polyadic tensor, $f(\theta)$, with $\theta \in \mathbb{S}^{n-1}$, see Eq. (1.4). It is interesting to introduce the multiplicities of its components, e.g. of $f_{\underbrace{i_1, i_1, i_2, i_2, i_1, i_3, i_1, i_2, i_5, i_4, \ldots}_{t \text{ indices}}}$, by introducing the $l_t$ integers $\{t(i_1), \ldots, t(i_{l_t})\}$, with $l_t \leqslant n$, so that

$$f_{i_1, i_2, \ldots, i_t}(\theta) = \prod_{q=1}^{l_t} \theta_{i_q}^{t(i_q)} \tag{2.3}$$

These $\{t(i_1), \ldots, t(i_{l_t})\}$, are by definition the multiplicities of index $i_1, \ldots, i_{l_t}$, for the component $(i_1, \ldots, i_t)$ of $f(\theta)$.
Note, that for any component, the sum of the multiplicities is $t$, the order of the tensor. Obviously, the knowledge of the index with their multiplicities entirely characterizes the corresponding component of the order $t$ tensor.

Using the notation presented Eq. (2.3) to denote the components of tensor $\theta' \circ \theta' \circ \cdots \circ \theta'$ and defining as multiplicities the $l$ numbers $\{s(i), s(j), \ldots\}$ [2], the components of $T'$ are :

$$T'_{ijiikjk\ldots} = \int_{\mathbb{S}^{n-1}} \theta'^{s(i)}_i \theta'^{s(j)}_j \theta'^{s(k)}_k \ldots g\left(\theta'_1\right) d\theta'. \tag{2.4}$$

---

[2] In the following, we shall denote $i, j, k, \ldots$ the component index, instead of $i_1, i_2, \ldots i_t$, so as to simplify the notation.

The components in both basis are related by:

$$T^{lmnop...} = R^l_i R^m_j R^n_k R^o_l R^p_m .... T'^{ijklm....}, \qquad (2.5)$$

where Einstein summation rule has been used, and $R^l_i$ denotes the matrix element $R(l,i)$.

The multiplicities $\{s(i), s(j), \ldots\}$ are the multiplicities of the index in the component of $T'$ in Eq. (2.4).

2.1. *Elementary results on $T'$*. In this section, we derive simple relations that will permit to simplify Eq. (2.5) and eventually obtain an analytical expression.
The multiplicity of indices strictly greater than one is an even number, otherwise the corresponding component of $T'$ is zero.

*proof.* Starting from the definition of $T'_{ijk...}$, Eq. (2.4), we first proceed to the integration over the index greater than one, $(\theta'_2, \ldots, \theta'_n)$. As $g$ is not involved, elementary considerations on spherical integration in $\mathbb{R}^n$ [5, Appendix D.] shows that the component $(i, j, \ldots)$ of $T'$ that have an odd exponent – among the list $(s(i), s(j), \ldots)$ – for at least one of their index higher than one, are zero.
On the contrary, components that have an even exponent – e.g. $\theta'^2_2, \theta'^6_3....$ – will give a non-zero value, upon integrating.
Note that $\theta'_1$ can have an odd or even multiplicity in T'.

Therefore, the parity of the tensor order, $t$ plays an important role. We shall distinguish both cases, make a detailed proof for $t$ even, and give the result for $t$ odd.

Let us consider $t$ is even, $t = 2p$. In this case, the multiplicity of $\theta'_1$ is even, and, from Lemma (2.1), a component of $T$, obtained from the non-zero components of $T'$ writes :

$$T^{lmnop....} = \sum_{i,j,k,l,m...=1}^{n} R^l_i R^m_j R^n_k R^o_l R^p_m ... T'^{ijklm....} = B_0 + B_2 + B_4 + ..., \qquad (2.6)$$

By definition, $B_a$ , with $a = (0, 2, \ldots, 2p)$ is the sum of all the terms of $T^{lmnop....}$ that contains $\theta'^a_1$ .

To separate the dependancy on the index of $T'$ and $R$, in Eq. (2.6), it is useful to notice that: All the coordinates in $T'$ are equivalent, up to some permutations.

Besides, the components have the properties:

- $T'^{ii...i} = T'^{22...2} \quad \forall i > 1$,

- $T'^{P(ii...ijj)} = T'^{22...233} \quad \forall i, j > 1$

- $T'^{P(ii...ijjkk)} = T'^{22...23344} \quad \forall i, j, k > 1$, where P is a permutation.

- ...

*proof.* From Eq. (2.4), a permutation in the coordinates of $T'$ do not change the value of the integrand. The value of $T'_{ijk...}$ do not depend on the order of the index $\{i, j, k, \ldots\}$

Then, using the change of variable defined using a rotation that brings $(\theta_i, \theta_j, \dots) \to (\theta_1, \theta_2, \dots)$, and using the invariance of the integration result over the sphere $\mathbb{S}^{n-1}$, one easily obtains the result.

Due to lemma (2.1), $B_0$ is the sum of terms of the type :

$$T'_{22\dots\dots 22},\; T'_{22\dots\dots 2233},\; T'_{22\dots 223344},\; \dots.$$

that can be factored out of the sum Eq. (2.6).

From Eq. (2.4), the components of $T'$ that appear in each of the $B_a$ can be expressed as

$$T'_{B_a} = I_{>1}\; G_a(t,n) \times S_{n-1} \Big/ \int_0^\pi \sin\phi_1^{n-2} d\phi_1 \tag{2.7}$$

$$G_a(t,n) = \int_0^\pi g(\cos\phi_1)\cos^a\phi_1 \sin\phi_1^{n-2+t-a} d\phi_1,$$

where $I_{>1}$ are the contribution of the indices strictly greater than one.

The expressions of $G_a(t,n)$ that permit, e.g. to compute the kernel that describes vectorial Electromagnetic Scattering [7] are given in Sec. (8).

Formally, one can write :

$$
\begin{aligned}
B_0 &= \sum_{i=2}^{n} \sum_{\substack{(\alpha,\beta,\dots,\varepsilon) \\ = P(i,i,\dots,i)}}^{n} \left(R_\alpha^l R_\beta^m R_\gamma^n R_\delta^o R_\varepsilon^p \dots\right) T'^{\alpha\beta\dots\varepsilon} \\
&+ \sum_{\substack{i,j=2 \\ i\neq j}}^{n} \sum_{\substack{(\alpha,\beta,\dots,\varepsilon) \\ = P(i,i,\dots,j,j)}}^{n} \left(R_\alpha^l R_\beta^m R_\gamma^n R_\delta^o R_\varepsilon^p \dots\right) T'^{\alpha\beta\dots\varepsilon} \\
&+ \dots
\end{aligned} \tag{2.9}
$$

Rearranging the term in the products, and using Lemma (2.1)

$$
\begin{aligned}
B_0 &= T'^{22\dots 22} \sum_{i=2}^{n} R_i^l R_i^m R_i^n R_i^o R_i^p \dots \\
&+ T'^{22\dots 2233} \sum_{\substack{i,j=2 \\ i\neq j}}^{n} \sum_{\substack{(\alpha,\beta,\dots,\varepsilon) \\ = P(l,m,n,o,p,\dots)}} R_i^\alpha R_i^\beta R_i^\gamma \dots R_j^\delta R_j^\varepsilon \\
&+ \dots\dots
\end{aligned} \tag{2.10}
$$

Note that $P$ does not cover all the possible permutations, but only those that do correspond to the terms that are actually present in the sum, Eq. (2.6). Indeed, among the $t!$ contributions coming from all the permutations on $(l,m,n,o,p,\dots)$, some should not be considered, as they would permute the upper index (of $R_a^b$) of matrix element that correspond to the same lower index. e.g., they would turn $R_i^l R_i^m R_i^n R_i^o R_j^p R_j^q$ to, e.g.

$R_i^n R_i^m R_i^l R_i^o R_j^q R_j^p$, not present in Eq. (2.6). This point should be taken care of, otherwise some contributions to $B_0$ would be counted several times.

In the general case, among the ensemble of the $t!$ permutations, $s(1)! s(2)! ... s(l)!$ of them correspond in fact to a single term in the sum Eq. (2.6).

Therefore, we define $P$ as the ensemble of permutations on $(l, m, n, o, p, ...)$, whose image is an ordered set, defined in the following way :

The image by P of the t-uplet $(l, m, n, ...)$ is a set $\{s(1); s(2); ...; s(l)\}$ made of $l$ multiplets, each of them being ordered, following the rule $l \leqslant m \leqslant n \leqslant o \leqslant p \leqslant q \leqslant ...$.
Note that $s(1) = 0$ for $B_0$.

Let us consider an example with $t = 6$, multiplicities $\{s(2) = 4; s(3) = 2\}$. The possible ordered permutations $P$ are defined by the following images of $(l, m, n, o, p, q)$ :

| $\{(l,m,n,o);(p,q)\}$ | $\{(l,n,o,p);(m,q)\}$ | $\{(l,o,p,q);(m,n)\}$ | $\{(m,n,o,p);(l,q)\}$ |
|---|---|---|---|
| $\{(l,m,n,p);(o,q)\}$ | $\{(l,n,o,q);(m,p)\}$ | | $\{(m,n,o,q);(l,p)\}$ |
| $\{(l,m,n,q);(o,p)\}$ | $\{(l,n,p,q);(m,o)\}$ | | $\{(m,n,p,q);(l,o)\}$ |
| $\{(l,m,o,p);(n,q)\}$ | | | $\{(m,o,p,q);(l,n)\}$ |
| $\{(l,m,o,q);(n,p)\}$ | | | |
| $\{(l,m,p,q);(n,o)\}$ | | | $\{(n,o,p,q);(l,m)\}$ |

In the general case, the number of term of each type in Eq. (2.6) is :

For term of the type $iiii...i$ for given $i$ : one term.

For term of the type $iiiijj$ for given $i$ and $j$ : : $\frac{(t)!}{(t-2)!2!}$

For term of the type $iiiijjkk$ for given $i, j, k$ : : $\frac{(t)!}{(t-4)!2!2!}$

Eventually, there are $C_1^{s(1)s(2)...s(l)} = \frac{(t)!}{s(1)!s(2)!...s(l)!}$ different permutations $P$ of the ensemble $(l, m, n, o, p, ...)$ that contribute to the sum of Eq. (2.6).

**3. Computation of $B_0$.** We will detail here the computation of $B_0$, and show that $B_n$, $n > 0$ can be expressed simply, as a function of $B_0$.

*Factorisation of the general expression.* For a given t-uplet among the $C_1^{s(1)s(2)...s(l)}$, we split it in $t/2$ doublets (let us recall that $t$ is even), and introduce all the circular permutations of these doublets. This will be helpfull to make appear the scalar product of two lines vectors of the rotation matrix $R$. For example, $\{(l, m, n, o); (p, q)\}$ is splitted in 3 doublets $\{(l, m); (n, o); (p, q)\}$, and we introduce the permutations corresponding to
- $\{(l, m, n, o); (p, q)\}$

- $\{(p, q, l, m); (n, o)\}$

- $\{(n,o,p,q);(l,m)\}$

In the general case, there are $C_2^{s(1)s(2)...s(l)} = \frac{(t/2)!}{(s(1)/2)!(s(2)/2)!...(s(l)/2)!}$ of such permutations.

Let us choose one t-uplet among the $C_1$ that are possible, and then define $a_i, b_i, c_i, \ldots$, $\forall i \in [2..n]$, taking pairwise :

$$\begin{aligned} a_i &= R_i^l R_i^m \\ b_i &= R_i^n R_i^o \\ c_i &= R_i^p R_i^q \\ &\ldots \end{aligned} \quad (3.1)$$

we can then rewrite the sum in Eq. (2.10) as

$$\begin{aligned} & \sum_{\substack{i,j=2 \\ i \neq j}}^{n} \sum_{\substack{(\alpha,\beta,\ldots,\varepsilon) \\ = \{(l,m,n,o);(p,q)\}}} R_i^\alpha R_i^\beta R_i^\gamma ... R_j^\delta R_j^\varepsilon \quad (3.2) \\ =\ & \sum_{\substack{i,j=2 \\ i \neq j}}^{n} \sum_{\substack{(\alpha,\beta,\ldots,\varepsilon) \\ = \{(l,m,n,o);(p,q)\}}} \frac{1}{C_2^{s(1)s(2)...s(l)}} (a_i b_i \ldots c_j + a_i b_j \ldots c_i + a_j b_i \ldots c_i) \end{aligned}$$

Therefore, we have in the general case

$$B_0 = T^{'22..22} \underbrace{\sum_{i=2}^{n} a_i b_i c_i ... e_i f_i}_{a\ single\ term} \quad (3.3)$$

$$+ \sum_{C_1\ terms} \frac{T^{'22..2233}}{C_2^{22..2233}} \sum_{i,j=2}^{n} (a_i b_i ... e_i f_j + a_i b_i ... e_j f_i + ... + a_j b_i ... e_i f_i) + ....$$

One notice that, for a given uplet among the $C_1$ possible t-uplet, one can construct different configurations, that would lead to different definitions of the $a_i$, $b_i$, ..., see Eq. (3.1), but would correspond to the same value for their products $a_i b_i c_j \ldots$.
For example, $\{(l,m,n,o)_i;(p,q)_j\}$ gives a contribution identical to $\{(m,o,l,n)_i;(p,q)_j\}$.
Let us establish the complete list of such equivalent t-uplets, splitted in pairs.
With our simple example, we obtain :

- $\{(l,m)_i\ (n,o)_i\ (p,q)_j\}$

- $\{(l,n)_i\ (m,o)_i\ (p,q)_j\}$

- $\{(l,o)_i\ (m,n)_i\ (p,q)_j\}$

- $\{(m,n)_i \ (l,o)_i \ (p,q)_j\}$

- $\{(m,o)_i \ (l,n)_i \ (p,q)_j\}$

- $\{(n,o)_i \ (l,m)_i \ (p,q)_j\}$

In the general case, we can count $C_3^{s(1)s(2)...s(l)}$ of such uplets that are equivalent to a single uplet of the initial list with $C_1$ elements. The total number of pairs for the sequence $\{lmnopq \ldots\}$ is

$$N = C_1.C_3 = \frac{(t)!}{\underbrace{(2)!\,(2)!...\,(2)!}_{t/2 \ times}} = \frac{(t)!}{2^{t/2}}$$

we can then simply rewrite

$$\sum_{C_1 \text{ terms}} \frac{T'^{22..2233}}{C_2^{22..2233}} \sum_{i,j=2}^{n} (a_i b_i ... e_i f_j + a_i b_i ... e_j f_i + ... + a_j b_i ... e_i f_i) + .... \quad (3.4)$$

$$= \sum_{C_1 \text{ terms}} \frac{1}{C_3^{22..2233}} \sum_{C_3 \text{ pairs}} \frac{T'^{22..2233}}{C_2^{22..2233}} \sum_{i,j=2}^{n} (a_i b_i ... e_i f_j + a_i b_i ... e_j f_i + ...)$$

Arranging the summations as

$$\sum_{C_1 \text{ terms}} \frac{1}{C_3^{22..2233}} \sum_{C_3 \text{ pairs}} \frac{T'^{22..2233}}{C_2^{22..2233}} \sum_{i,j=2}^{n} (...) \quad (3.5)$$

$$= \frac{1}{C_3^{22..2233}C_1^{22..2233}} \sum_{C_1 \text{ terms}} \sum_{C_3 \text{ pairs}} \frac{T'^{22..2233}C_1^{22..2233}}{C_2^{22..2233}} \sum_{i,j=2}^{n} (....)$$

One obtains

$$B_0 = \frac{1}{N} \sum_{\substack{\text{all the pairs} \\ \text{in } (l,m,n,o,p,...)}} \left[ T'^{22..22} \underbrace{\sum_{i=2}^{n} a_i b_i c_i ... e_i f_i}_{a \ single \ term} \right. \quad (3.6)$$

$$+ \left. \frac{T'^{22..2233}}{C_2^{22..2233}} C_1^{22..2233} \sum_{i,j=2}^{n} (a_i b_i ... e_i f_j + a_i b_i ... e_j f_i + ... + a_j b_i ... e_i f_i) + ... \right]$$

The rationale behind this rearrangement of the sum is to proceed to a factorization of the type :

$$\sum_{i=2}^{n} a_i b_i c_i ... e_i f_i + \sum_{i,j=2}^{n} (a_i b_i ... e_i f_j + a_i b_i ... e_j f_i + ... + a_j b_i ... e_i f_i) + ...$$

$$= \left(\sum_{i=2}^{n} a_i\right)\left(\sum_{i=2}^{n} b_i\right)...\left(\sum_{i=2}^{n} f_i\right) \tag{3.7}$$

To do so, one needs the following theorem to be proved : any component of the tensor $\frac{T'C_1}{C_2}$ only depends on $n$ and $t$, but not on the values of the $s(i)$ in the configuration $\{s(1), s(2), \ldots, s(l)\}$.

*Proof.* As $C_1$ and $C_2$ are known, we need to compute T' in the general case.

A component of $T'$ of order $t$, defined by its values $s(i)$ with $1 < i < l$ can be written :

For $l > 2$ :
$$T' = G_0(t,n) \times S_{n-1} / \int_0^\pi \sin \phi_1^{n-2} d\phi_1 \times D \times E$$

For $l = 2$ :
$$T' = G_0(t,n) \times S_{n-1} / \int_0^\pi \sin \phi_1^{n-2} d\phi_1 \times E$$

And for $l = 1$ :
$$T' = G_0(t,n) \times S_{n-1} / \int_0^\pi \sin \phi_1^{n-2} d\phi_1$$

With
$$D = \prod_{j=2}^{l-1} \left[\int_0^\pi \cos \phi_j^{s(j)} \sin \phi_j^{s(j+1)+....+s(l)} \sin \phi_j^{n-j-1} d\phi_j / \int_0^\pi \sin \phi_j^{n-j-1} d\phi_j\right]$$

$$E = \int_0^\pi \cos \phi_l^{s(l)} \sin \phi_l^{n-l-1} d\phi_l / \int_0^\pi \sin \phi_l^{n-l-1} d\phi_l$$

Let us consider the case $(l > 2)$ :

Setting $u(j) = \sum_{l=j+1}^{l} s(l)$, one gets

$$D = \prod_{j=2}^{l-1} \left(\frac{s(j)-1}{n-j-1+u(j-1)} \times \frac{s(j)-3}{n-j-1+u(j-1)-2} \cdots \frac{1}{n-j-1+u(j)+2}\right) \times$$
$$\left(\frac{n-j-1+u(j)-1}{n-j-1+u(j)} \times \frac{n-j-1+u(j)-3}{n-j-1+u(j)-2} \cdots \frac{n-j-1+1}{n-j-1+2}\right)$$
$$= \prod_{j=2}^{l-1} \frac{(s(j)-1)!!}{(n-j-1+u(j-1))!!} \frac{(n-j-1+u(j)-1)!!}{(n-j-2)!!} (n-j-1)!!$$

One can show easily that :

$$\prod_{j=2}^{l-1} \frac{(n-j-1)!!}{(n-j-2)!!} = \frac{(n-3)!!}{(n-1-l)!!}$$

and
$\prod_{j=2}^{l-1} \frac{(n-j-1+u(j)-1)!!}{(n-j-1+u(j-1))!!} = \frac{(n-l-1+u(l-1))!!}{(n-3+u(1))!!}$ with $u(1) = t$ and $u(l-1) = s(l)$

In the same way,

$$E = \frac{(s(l)-1)!!}{(n-l-1+s(l))!!}(n-l-1)!!$$

So that,

$$D \times E = \frac{(n-3)!!}{(n-3+t)!!} \prod_{j=2}^{l} (s(j)-1)!!$$

And therefore,

$$T' = G_0(t,n) \times \left(S_{n-1}/\int_0^\pi \sin\phi_1^{n-2} d\phi_1\right) \frac{(n-3)!!}{(n-3+t)!!} \prod_{j=2}^{l} (s(j)-1)!!$$

As seen here, the components of $T'$ depend on the $\{s(i)\}$, so that the factorization of $B_0$ is not straightforward. We will now investigate $T'\frac{C_1}{C_2}$, using that

$$C_1/C_2 = \left(\frac{(t)!}{\prod_{j=2}^{l} s(j)!}\right) / \left(\frac{(t/2)!}{\prod_{j=2}^{l} (s(j)/2)!}\right)$$

Note that[3], if $l = 2$, $C_1/C_2 = 1$.

To do so, we first focus on the function

$$f(n) = 2^{n/2}(n-1)!!/n!$$

Obviously,

$$f(n+2) = \frac{2}{n+2}f(n)$$

So that

$$f(n) = 1/(n/2)!$$

Consequently,

$$C_1/C_2 = \left(\frac{(t)!}{\prod_{j=2}^{l} s(j)!}\right) / \left(\frac{(t/2)!}{\prod_{j=2}^{l} (s(j)/2)!}\right) = \frac{2^{t/2}(t-1)!!}{\prod_{j=2}^{l} 2^{s(j)/2}(s(j)-1)!!} = \frac{(t-1)!!}{\prod_{j=2}^{l} (s(j)-1)!!}$$

and

$$T'C_1/C_2 = G_0(t,n) \times \left(S_{n-1}/\int_0^\pi \sin\phi_1^{n-2} d\phi_1\right) \frac{(n-3)!!}{(n-3+t)!!}(t-1)!!$$

We have just shown that $\frac{T'C_1}{C_2}$ only depends on $n$ and $t$, but not on the values of the $s(i)$ in the configuration $\{s(1), s(2), \ldots, s(l)\}$. This relation is also true for $l = 2$.

---

[3]We recall that $s(1) = 0$, for $B_0$.

This proves theorem (3).

We also note that :

$$T'C_1/C_2 = T'^{22334455....}C_1^{22334455...}/C_2^{22334455...} = T'^{22334455....}(t)!/2^{t/2}/(t/2)!$$

Finally,

$$B_0(t,n) = \frac{T'^{223344....}}{(t/2)!} \sum_{\substack{all\,the\,pairs \\ a,b,c,....}} \left(\sum_{i=2}^n a_i\right)\left(\sum_{i=2}^n b_i\right).....\left(\sum_{i=2}^n f_i\right)$$

What can be written explicitely,

$$B_0(t,n) = G_0(t,n)\left(S_{n-1}/\int_0^\pi \sin\phi_1{}^{n-2}d\phi_1\right)\frac{(n-3)!!}{(n-3+t)!!}\frac{2^{t/2}}{t!!} \qquad (3.8)$$

$$\times \sum_{\substack{all\,the\,pairs \\ a,b,c,....}} \left(\sum_{i=2}^N a_i\right)\left(\sum_{i=2}^N b_i\right).....\left(\sum_{i=2}^N f_i\right)$$

For the pair (p,q), using the properties of rotation matrix Sec. (7),

$$\sum_{i=2}^n a_i^{pq} = \sum_{i=2}^n R_i^p R_i^q = -\xi_p\xi_q + \delta_{pq}$$

Now that $B_0$ has been calculated, we will see that we can calculate simply all the $B_n, n > 0$.

*Generalization to $B_n$.* Getting back to Eq. (2.6), one note that $B_2$ will write :

$$B_2 = \sum_{(\alpha,\beta)\in(l,m,n,o,p,...)} R_1^\alpha R_1^\beta \, J_2(t,n)|_{(\gamma,\delta,\varepsilon...\eta,\mu)\neq(\alpha,\beta)} \qquad (3.9)$$

with

$$\begin{aligned} J_2 &= T'^{1122..22}\sum_{i=2}^n R_i^l R_i^m R_i^n \ldots \\ &+ T'^{1122..2233}\sum_{\substack{i,j=2 \\ j>i}}\sum_{\substack{(\alpha,\beta,\ldots,\varepsilon) \\ =P(l,m,n,o,p,\ldots)}} R_i^\gamma R_i^\delta R_i^\varepsilon...R_j^\eta R_j^\mu \ldots \\ &+ \ldots \end{aligned} \qquad (3.10)$$

As demonstrated previously for $B_0$, one readily sees that $J_2$ only depends on $n$ and $t$ but not on the multiplicities $\{s(i)\}$. This uses the fact, see Eq. (2.7), that $T'^{1122334...}$ is separable, in a product of the type $A^{11}B^{22334...}$

Considering now $B_i$, $i$ even, one notice that many simplifications can be done. Firstly, the sum on the products of $R_a^b$ are of the same type for all the $B_i$. Secondly, only the number of term is different : $t = 2p$ terms for $B_0$, $t - 2$ for $B_2$, ... .

From all this, we deduce that

$$J_a\left(t,n\right) = \frac{G_a\left(t,n\right)}{G_0\left(t-a,n\right)} J_0\left(t-a,n\right) \qquad a \leqslant 2p, \text{ a } is\ even, \tag{3.11}$$

where $J_0\left(t,n\right) = B_0$.

3.1. *Final Results.* Finally, let us start from,

$$\begin{aligned} T^{lmnop....} &= \sum_{i,j,k,l,m...=1}^{n} R_i^l R_j^m R_k^n R_l^o R_m^p ... T^{'ijklm....} \tag{3.12} \\ &= J_0\left(t,n\right) + \sum_{(\alpha,\beta)\in(l,m,n,o,p,...)} \xi_\alpha \xi_\beta \left. J_2\left(t,n\right)\right|_{(\gamma,\delta,\varepsilon...\eta,\mu)\neq(\alpha,\beta)} + ..., \end{aligned}$$

and denote $P_{2w}$ the ensemble of the $2w$-uplets extracted from the $t$-uplet (lmnop....), and $\overline{P_{2w}}$ the complementary of $P_{2w}$ in $P_t$, one can write, for even $t$ :

$$\begin{aligned} T^{lmnop....} &= \sum_{w=0}^{t/2} \sum_{P_{2w}} \underbrace{\xi_\alpha \xi_\beta \cdots}_{2w \text{ terms},(\alpha,\beta..)\in P_{2w}} J_{2w,\overline{P_{2w}}}\left(t,n\right) \tag{3.13} \\ &= \sum_{w=0}^{t/2} \sum_{P_{2w}} \underbrace{\xi_\alpha \xi_\beta \cdots}_{2w \text{ terms},(\alpha,\beta..)\in P_{2w}} \frac{G_{2w}\left(t,n\right)}{G_0\left(t-2w,n\right)} J_{0,\overline{P_{2w}}}\left(t-2w,n\right), \tag{3.14} \end{aligned}$$

where $J_{0,\overline{P_{2w}}}$ is computed using Eq. (3.8), considering only the pairs of $\overline{P_{2w}}$ to define the $a_i, b_i, c_i, \ldots$ in Eq. (3.1).

The extension of this work to odd $t$ is direct and easy. We obtain :

$$T^{lmnop....} = \sum_{w=0}^{(t-1)/2} \sum_{P_{2w+1}} \underbrace{\xi_\alpha \xi_\beta \cdots}_{(2w+1) \text{ terms},(\alpha,\beta..)\in P_{2w+1}} J_{2w+1,\overline{P_{2w+1}}}(t,n) \qquad (3.15)$$

$$= \sum_{w=0}^{(t-1)/2} \sum_{P_{2w+1}} \underbrace{\xi_\alpha \xi_\beta \cdots}_{(2w+1) \text{ terms},(\alpha,\beta..)\in P_{2w+1}} \frac{G_{2w+1}(t,n)}{G_0(t-2w-1,n)} J_{0,\overline{P_{2w+1}}}(t-2w-1,n)$$

$$(3.16)$$

One can therefore conclude that the FT of a kernel defined by Eq. (1.3) with a polyadic tensor $f$ is a polyadic tensor, whose coefficients can be computed by Eqs (3.14, 3.16).

**4. Comparaison with numerical simulations.** In this section, we compare the direct evaluation of $T$, by numerical integration, see Eq. (1.7), and using our expression Eq. (3.14, 3.16). Note that we limit ourselves to odd values of $t$ in this section. The extension to even values is straightforward. Besides, a recursive program, that works for odd and even values of $t$ is available [8]. It permits to compute the coefficients in the sum Eq. (3.14, 3.16).

For odd values of $t$, we shall work with the function : $g(x) = sgn(x)$.
The goal is to compute any component of the tensor $T$ of order $t$, for a given unitary vector $\xi \in \mathbb{R}^n$.

4.1. *Analytic calculation based on recursivity.*

4.1.1. *odd values of $t$.* The tensor $T$ is a sum, whose general structure, *for odd values of $t$*, can be described as :

$$T(n,t,w,E,E_c) = S_{n-1} Z(t,n,w) \underbrace{\xi_\alpha \xi_\beta \cdots}_{2w+1 \text{ terms},(\alpha,\beta..)\in P_{2w+1}} \underbrace{(-\xi_p\xi_q + \delta_{pq})(-\xi_r\xi_t + \delta_{rs}) \cdots}_{(p,q,s,t,...)\in \overline{P_{2w+1}}}$$

with $E \in P_{2w+1}$, $E_c \in \overline{P_{2w+1}}$ and

$$Z(t,n,w) = \frac{(n-3)!! (2w)!!}{(n-2+t)!! (t-2w-1)!!} \frac{2^{(t-2w-1)/2}}{\int_0^{\pi/2} \sin\phi_1^{n-2} d\phi_1},$$

Where we denote $S_{n-1} = \frac{2\pi^{n/2}}{\Gamma(n/2)}$.
A product of factors $\xi_\alpha$ and $(-\xi_p\xi_q + \delta_{pq})$ appears. Therefore, one can use a recursive algorithm :

Starting from the term $w = (t-1)/2$ corresponding to the "level zero" :
$S_{n-1} Z(t,n,(t-1)/2) \underbrace{\xi_\alpha \xi_\beta \cdots}_{t \text{ terms},(\alpha,\beta..)\in P_t}$, the general term of level 1 is constructed, taking a couple of $(\xi_p, \xi_q)$ among the $t$ possible $\xi_j$, and replacing the product $\xi_p\xi_q$ by $(-\xi_p\xi_q + \delta_{pq})$. This procedure is repeated up to terms that have only one factor $\xi$, corresponding to the level $(t-1)/2$. Note that the arborescent structure corresponding to the double summation $\sum_{w=0}^{(t-1)/2} \sum_{P_{2w+1}}$ should be fully explored.

$T(n,t,w,E,E_c)$ is computed at each node of the graph. Finally, $Z$ depends on $w$ and should be modified at each recursion level. In the implementation, all the needed values of $Z(t,n,w)$ are computed at the begining, before calling the recursive routine.

The sought value of $T$ component will be the sum of all the intermediate values computed at each node of the graph. Finally, let us mention that to each component correspond in general a different arborescent structure.

4.1.2. *even values of t.* For even values of $t$, the procedure is similar, with the definitions :

$$T(n,t,w,E,E_c) = S_{n-1} Z(t,n,w) \underbrace{\xi_\alpha \xi_\beta \ldots}_{2w \text{ terms}, (\alpha,\beta..) \in P_{2w}} \underbrace{(-\xi_p \xi_q + \delta_{pq})(-\xi_r \xi_t + \delta_{rs}) \ldots}_{(p,q,s,t,\ldots) \in \overline{P_{2w}}}$$

with $E \in P_{2w}$, $E_c \in \overline{P_{2w}}$ and

$$Z(t,n,w) = \frac{(n-3)!! 2^{(t-2w-4)/2}}{(n-3+t-2w)!!(t-2w)!!} \frac{\Gamma\left(\frac{2w+1}{2}\right)\Gamma\left(\frac{n+s-(2w+1)}{2}\right)}{\Gamma\left(\frac{n+s}{2}\right)} \times \\ \left(\Psi\left(\frac{n+s}{2}\right) - \Psi\left(\frac{2w+1}{2}\right)\right) \frac{1}{\int_0^{\pi/2} \sin\phi_1^{n-2} d\phi_1}, \quad (4.1)$$

where $\Psi(x) := \frac{\Gamma'(x)}{\Gamma(x)}$ is the digamma function.

4.2. *Numerical simulation.* In the following, the results are expressed as the ratio of the $T$ component by the unitary sphere surface $S_{n-1}$.

To make the direct integration of Eq. (1.7), we use Monte-Carlo methods.

In a first approach (labeled MC1, hereafter), we consider a sphere of unitary radius in $\mathbb{R}^n$, and generate a large number $N$ of points, uniformly spread on the sphere surface [9]. To do so, we use a Müller algorithm [10]. At each point, the components of tensor $f$ are calculated. The integral then evaluated as :

$$T(\xi)/S_{n-1} = \frac{1}{S_{n-1}} \int_{\mathbb{S}^{n-1}} f(\theta) g(\xi \bullet \theta) d\theta \approx \frac{1}{N} \sum_{i=1}^{N} f(\theta^i) g(\xi \bullet \theta^i)$$

We have also used two other Monte-Carlo methods – labeled MC2 and MC3 –, in which the points are spread differently over the sphere surface :

MC2 : Instead of the Müller algorithm, we draw the spherical coordinates of the points on the sphere: $\theta \in [0..\pi]$ and $\phi \in [0..2\pi]$, using pseudo-random numbers – *rand* function of matlab.

MC3 : The spherical coordinates $(\theta, \phi)$ are drawn as quasi-random numbers, using matlab function *haltonset* [11], that generates a Halton set of numbers more evenly distributed than the pseudo random numbers.

4.3. *Results.* For example, we compare the numeric and exact results for a non trivial case, where $n = 3$ and $t = 5$. We used typically $N = 5 \times 10^6$ points on the sphere to compute the integrals. Three components of the tensor have been computed, see Figures (1-3). In most cases, one notice at first sight a good agreement, as the relative discrepancy is less than $2 \times 10^{-3}$. This is a check of the obtained formulas Eq. (3.16).

However, a carefull study shows that the direct numeric evaluation through Monte-Carlo simulation is difficult and unaccurate in some cases.

4.3.1. *accuracy of the numerical solution.* As expected, we note that the difference between exact and numerical results varies as $N^{-\frac{1}{2}}$.

For a given $N$, one observes a relative discrepancy between exact and numerically integrated result, which is greater when the computed component has a small value, whereas the absolute discrepancy is roughly uniform vs $\xi$, see fig.(1b, 2b, 3b), and comes from the Monte-Carlo method itself. A systematic study versus $N$ confirms this observation.

Tab.(1) shows an example where the numerical integration is very unaccurate. The absolute value of $T^{13333}/S_2$ is computed using MC1, MC2, and MC3, close to the pole $\xi_1 = 0$ of the sphere. One readily sees that even with $N = 5 \times 10^7$ points, the accuracy is not satisfying, even using the Halton method. We have examined the accuracy for other locations on the unit sphere. The results Fig. (4) show that even with $N = 5 \times 10^7$ points, and using MC3 method, a non negligible fraction of the directions $\xi$ are badly computed. Finally, let us note that the execution time of the recursive program that computes $T^{lmno\cdots}/S_{n-1}$ is order of magnitudes faster than the direct integration. Typically, when one needs 0.2 s to make the analytic computation, one needs 200 s with a MC3 and $5 \times 10^7$ points to carry out the numerical integration.

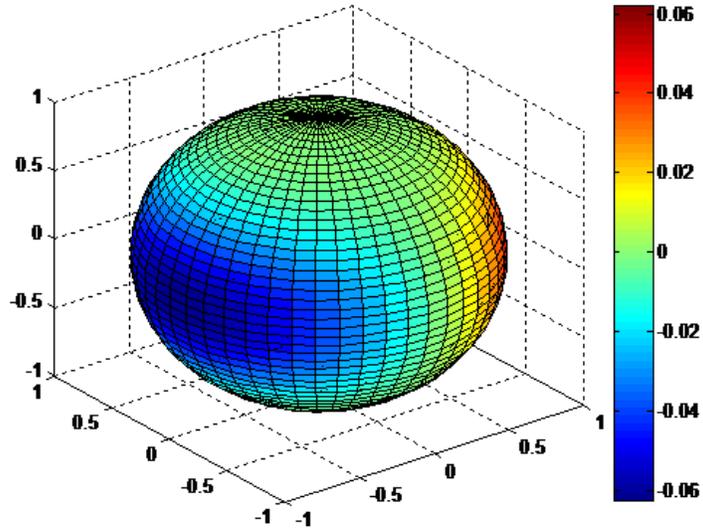

(A) Display of the exact result, Eq. (3.16)

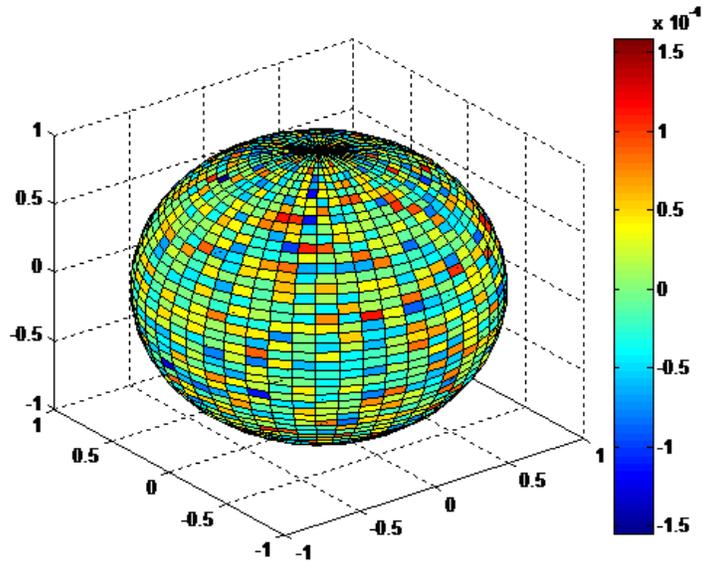

(B) Absolute difference between exact calculation and Monte Carlo (MC3) simulation, corresponding to fig.(1a).

FIG. 1. Comparison of exact and simulated result for $T^{13333}/S_2$, corresponding to the tensor $f_{13333}(\theta) = \theta_1 \theta_3^4$ as a function of unitary vector $\xi$ on a $\mathbb{R}^3$ sphere.

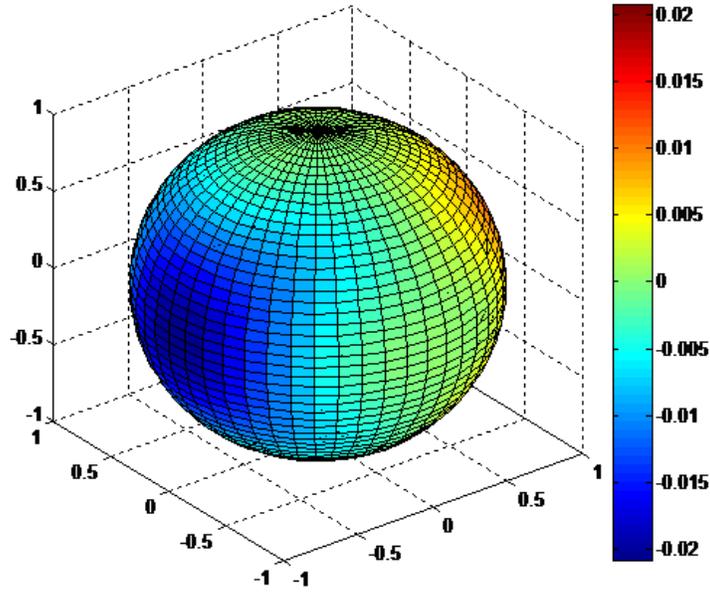

(A) Display of the exact result, Eq. (3.16)

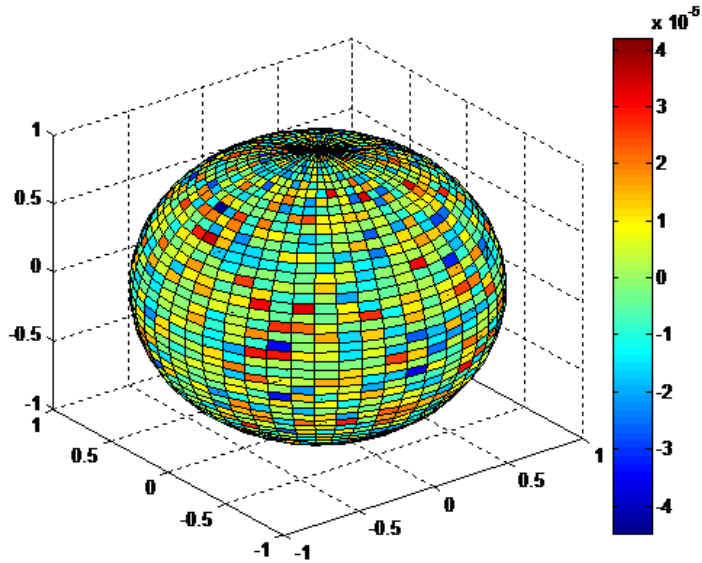

(B) Absolute difference between exact calculation and Monte Carlo (MC3) simulation, corresponding to fig.(2a).

FIG. 2. Comparison of exact and simulated result for $T^{12233}/S_2$, corresponding to the tensor $f_{12233}(\theta) = \theta_1 \theta_2^2 \theta_3^2$ as a function of unitary vector $\xi$ on a $\mathbb{R}^3$ sphere.

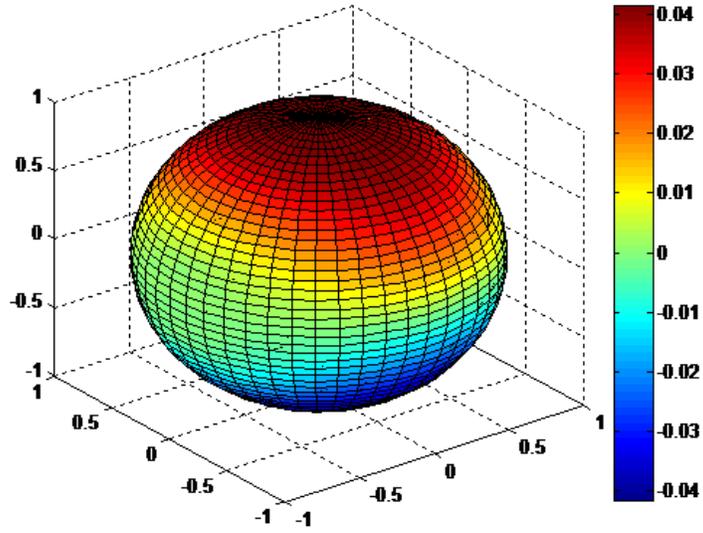

(A) Display of the exact result, Eq. (3.16)

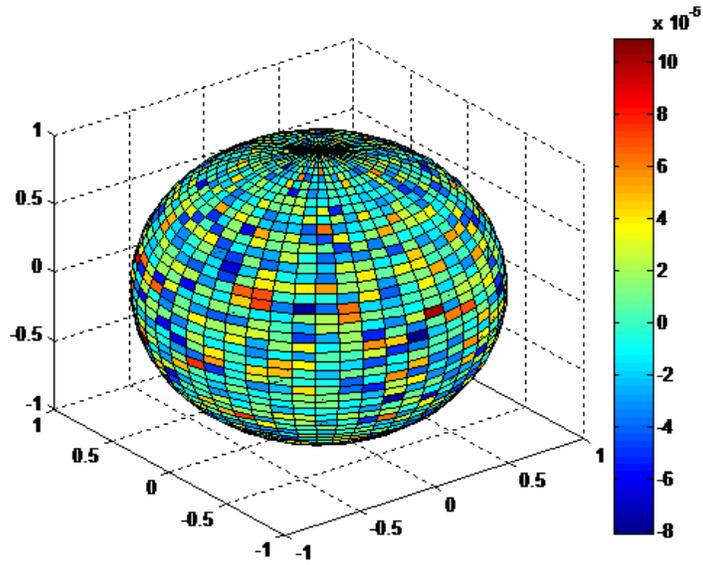

(B) Absolute difference between exact calculation and Monte Carlo (MC3) simulation, corresponding to fig.(5).

FIG. 3. Comparison of exact and simulated result for $T^{11333}/S_2$, corresponding to the tensor $f_{11333}(\theta) = \theta_1^2 \theta_3^3$ as a function of unitary vector $\xi$ on a $\mathbb{R}^3$ sphere.

| Method/N | 50 000 | 500 000 | 5 000 000 | 50 000 000 |
|---|---|---|---|---|
| MC1 | 4 $10^{-4}$ | 5 $10^{-5}$ | 3 $10^{-5}$ | 1 $10^{-5}$ |
| MC2 | 3 $10^{-4}$ | $10^{-4}$ | 6 $10^{-5}$ | 5 $10^{-6}$ |
| MC3 | 3 $10^{-6}$ | $10^{-6}$ | -1.64 $10^{-7}$ | -1.56 $10^{-7}$ |

Tab.(1) Absolute difference between exact calculation and Monte-Carlo simulation, for the different Monte-Carlo method MC1, MC2, MC3, at the point $\xi = (-0.0054, 0.1491, 0.9888)$. The exact value is $T^{13333}/S_2 = -1.67\ 10^{-7}$.

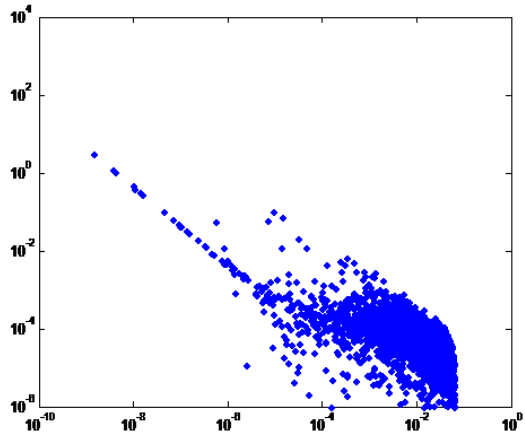

FIG. 4. Relative difference between exact calculation and Monte Carlo MC3 simulation, for $T^{13333}/S_2$, plotted versus the exact value of $T^{13333}/S_2$, for different points on the sphere surface.

**5. Discussion.** Let us give an exemple, in the field of image processing. We therefore take $n = 2$, and consider an image, represented by a function $u \in \mathcal{S}(\mathbb{R}^2)$. We shall consider the following kernel:

$$f(\theta) = \theta_1^{t(1)} \theta_2^{t(2)} = \cos(\phi)^{t(1)} \sin(\phi)^{t(2)}, \qquad (5.1)$$

in the case where $t = t(1) + t(2)$ is even. There, $\phi$ is the polar angle in the plane, by respect to some reference direction (here, the horizontal axis $x_1$). The goal is to compute how the 2D Riesz transformation defined by the convolution of $u$ with the distribution Eq. (1.3) will modify the image.

From Eq. (5.1), one has $\int f(\theta) d\theta = 0$, so that hypothesis Eq. (1.5) is fullfiled. With our method, we can rapidly compute the Fourier transform of the kernel and obtain the transformed image:

$$\tilde{u} = \mathcal{F}^{-1}\left[\widehat{W_f}.\hat{u}\right], \qquad (5.2)$$

where $\mathcal{F}^{-1}$ is the inverse Fourier transform. The analytic expression of $\widehat{W_f}$ obtained using Eq. (3.14) can be expressed more explicitely as:

$$\widehat{W_f}(\xi) = \sum_{p=0}^{t(1)} \sum_{q=0}^{t(2)} \left(\frac{\xi_1}{\|\xi_1\|}\right)^{p+t(2)-q} \left(\frac{\xi_2}{\|\xi_2\|}\right)^{q+t(1)-p} C_{t(1)}^p C_{t(2)}^q B_{p,q}, \quad (5.3)$$

where

$$B_{p,q} = (-1)^{[t(1)-p]} \int_0^{2\pi} [\cos(\alpha)]^{p+q} [\sin(\alpha)]^{t-p-q} \ g(\cos(\alpha)) \ d\alpha \quad (5.4)$$

Introducing the angle $\nu \in [0..2\pi]$ that describes the orientation of the vector $\xi \in \mathbb{R}^2$. One can define a family of rotated kernels $\widehat{W_{f,\theta_0}}, \theta_0 \in [0..2\pi]$ by:

$$\widehat{W_{f,\theta_0}}(\nu) = \sum_{p=0}^{t(1)} \sum_{q=0}^{t(2)} (\cos(\nu-\theta_0))^{p+t(2)-q} (\sin(\nu-\theta_0))^{q+t(1)-p} C_{t(1)}^p C_{t(2)}^q B_{p,q}. \quad (5.5)$$

Using the kernel Eq. (5.5), we eventually show Fig. (5), one can locate easily the corners in a picture made of overlapping rectangles. In order to do so, we choose the specific kernel defined by $t(1) = 3$, $t(2) = 1$.

This method we propose avoids, for example, to use the more complex technique of the monogenic image, which necessitate to make heavier computations [12] to locate intrinsically 2D features (such as corners).

Multimensional Riesz transforms for image processing have been well studied, e.g. in the community of computer vision. However, they have always been defined (to our knowledge) directly in the reciprocal space, as a Fourier multiplier [13]. The present work also permits to make the link between the singular integral transform in direct space, in which the image is defined, Eq. (1.3), and these Fourier mutlipliers. This relationship might help, e.g. to compare the different numerical methods for image processing, as, currently, some are constructed as convolutions in the real space (Harris corner detector [14], Sobel filters, ...) and others as multiplication in the reciprocal space.

Starting from this example, the extension to $\mathbb{R}^3$, with application to e.g. tomography is straightforward.

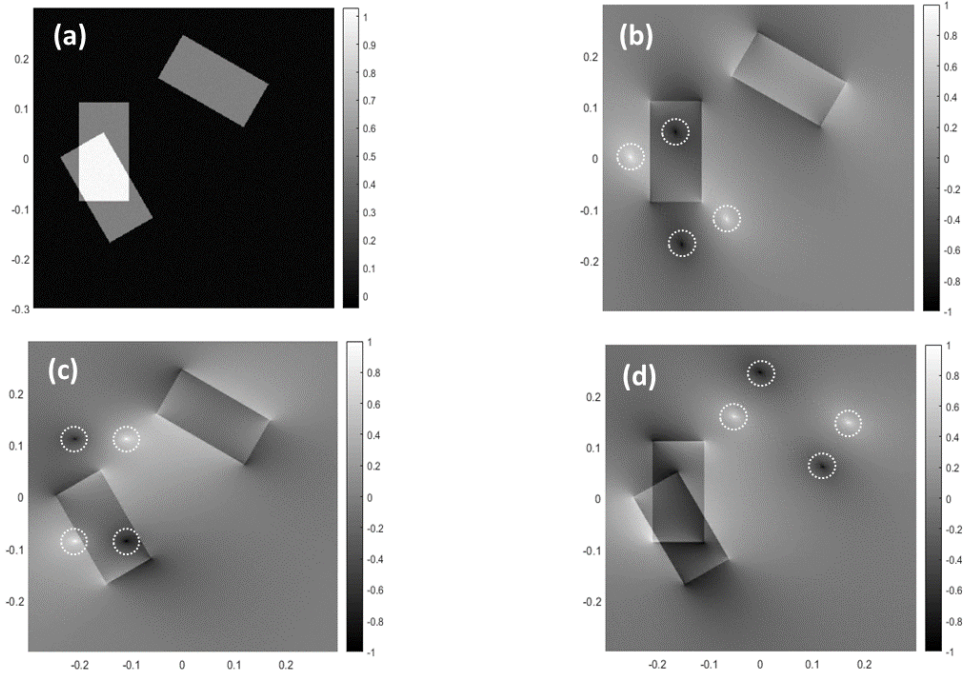

Fig. 5. Multidimensional Riesz transform of a simple image. The singular kernel Eq. (5.1), with $t(1) = 3$, $t(2) = 1$ is used in the distribution Eq. (1.3), to operate on the original image displayed panel (a). Panel (b-d) display the modified images obtained by convolution using Eq. (5.2), respectively for $\theta_0 = \pi/3$, $\theta_0 = \pi/2$, $\theta_0 = \pi/6$. The rectangle corners appear as bright or dark spots in the transformed images, enabling the selective detection of the corners of inclined rectangles (dotted circles have been added as a guide to the eyes).

**6. Conclusion.** In summary, we have generalized the Riesz transform to polyadic tensorial integral transforms of higher order. This somehow completes the work of Adkins [15] on integrals with polyadic kernels that are weakly singular. For all tensor orders, we provide an analytic expression, tested in simple cases against numeric calculations. This enables fast computation of all the components of the kernel Fourier transform. This permits to define the transform as a convolution, under the hypothesis of $\int_{\mathbb{S}^{n-1}} f = 0$.
As a perspective, we could mention the possible use of this integral transform for data analysis. For two-dimensional problems (e.g. image analysis), we have discussed how we could spot 2D features in an image (corners). It would be interesting in the future to compare this transform defined in the real space, to the Riesz transform defined in the reciprocal space, that is usually used in the computer vision community. Generalization to fractional singular transforms [16], with the goal of improving image analysis in the presence of noise is also a possibe route to follow. Beyond two-dimensional problems, let us note that, in association with a Canonical Polyadic decomposition [6,17,18], our work could permit to handle analytically the convolution with singular kernels of much larger dimension, so as to look for patterns in "big datas" arrays.

**Acknowledgements.** We thank Johan Debayle for helpful discussions.

## 7. Annex 1 : change of base matrix.

We operate a change of basis, so that the new basis has $\vec{\xi}$ for first (unitary) vector.

This transformation is :
$$\left(\vec{i_1}, \vec{i_2}, .....\vec{i_k},....\vec{i_n}\right) \to \left(\vec{e_1} = \vec{\xi}, \vec{e_2}, .....\vec{e_k},....\vec{e_n}\right)$$

the former basis being orthonormal.

The new basis is constructed by Schmidt orthonormalization :
$$\vec{e_1} = \vec{\xi}$$
$$\vec{e_2} = \vec{i_1} - \left(\vec{i_1} \bullet \vec{\xi}\right) \vec{\xi}$$
$$.$$
$$.$$
$$\vec{e_k} = \vec{i_{k-1}} - \sum_{i=1}^{k-1} \frac{\left(\vec{i_{k-1}} \bullet \vec{e_i}\right)}{\vec{e_i} \bullet \vec{e_i}} \vec{e_i}$$
$$.$$
$$.$$

At this stage, the $\vec{e_k}$ are not unitary.

Introducing $B(j) = \sqrt{1 - \sum_{k=1}^{j} \xi_k^2}$ with $B(0) = 1$, an easy computation gives that :
$$\|\vec{e_k}\| = B(k-1)/B(k-2)$$

The change of basis matrix is defined as $[\vec{e_1}, \vec{e_2}/\|\vec{e_2}\|, .....\vec{e_k}/\|\vec{e_k}\|, ....\vec{e_n}/\|\vec{e_n}\|]$, where the $\vec{e_k}/\|\vec{e_k}\|$ are unitary columns vectors.

Let us now denote $\vec{e_k}$ the normalized vectors. The components of $\vec{e_k}$ expressed in the former basis are :
$$\vec{e_1} = \vec{\xi}$$
$$\vec{e_k} = \begin{pmatrix} 0, 0, ...e_k^{k-1} = B(k-1)/B(k-2), e_k^k = -\xi_{k-1}\xi_k/(B(k-1)B(k-2)), \\ ..., e_k^j = -\xi_{k-1}\xi_j/(B(k-1)B(k-2)), .... \end{pmatrix} k>1$$

The matrix $R$ is therefore triangular. It has the property that :

If $\vec{e_k}$ is the k$^{\text{th}}$ column vector of $R$, and if $\vec{E_j}$ is the j$^{\text{th}}$ line vector of $R$, then
$$\vec{E_j} \cdot \vec{e_k} = \delta_{j,k}.$$

## 8. Annex 2 : Computation of $G_a(t, n)$.

We first recall the definition:
$$G_a(t, n) = \int_0^\pi g(\cos \phi_1) \cos^a \phi_1 \sin^{n-2+t-a} \phi_1 \, d\phi_1$$

Let us consider $g(x) = sgn(x)$

We get :
$$G_a(t, n) = \begin{cases} 2 \int_0^{\pi/2} \cos^a \phi_1 \sin^{n-2+t-a} \phi_1 \, d\phi_1 & a \text{ odd} \\ 0 & a \text{ even} \end{cases}$$

using the relation [19], for even $\beta$ :
$$\int_0^{\pi/2} (\sin \theta)^\alpha (\cos \theta)^\beta \cos \theta \, d\theta = \frac{\beta(\beta-2)......2}{(\alpha+\beta-1)...(\alpha+3)(\alpha+1)} \frac{1}{(\alpha+\beta+1)}$$

24 M. PERRIN AND F. GRUYOne gets :

$$G_a(t,n) = 2\frac{(a-1)(a-3)...2}{(n-2+t)(n-2+t-2)...(n-2+t-a+1)} = 2\frac{(a-1)!!(n-2+t-a-1)!!}{(n-2+t)!!},$$

what is also valid for $a = 1$.

Considering now $g(x) = -\ln(|x|)$, we obtain :

$$G_a(t,n) = \begin{cases} -2\int_0^{\pi/2} \ln(|\cos\phi_1|) \cos^a \phi_1 \sin \phi_1^{n-2+t-a} d\phi_1 & a \text{ even} \\ 0 & a \text{ odd} \end{cases}$$

So that

$$G_a(t,n) = \frac{1}{2}\frac{\Gamma\left(\frac{a+1}{2}\right)\Gamma\left(\frac{n+t-(a+1)}{2}\right)}{\Gamma\left(\frac{n+t}{2}\right)}\left(\Psi\left(\frac{n+t}{2}\right) - \Psi\left(\frac{a+1}{2}\right)\right)$$

for $a$ even, where $\Psi(x) := \frac{\Gamma'(x)}{\Gamma(x)}$ is the digamma function.